\documentclass[12pt]{article}
\usepackage{amscd,amsfonts,amsmath,amssymb,yfonts}
\usepackage{amstext,amsthm,mathbbol,mathrsfs,stmaryrd} 
\usepackage[cmtip,arrow]{xy}
\usepackage{pb-diagram,pb-xy} 
\usepackage{curves,graphicx,texdraw} 
\usepackage{tikz} 
\usetikzlibrary{calc} 
\usepackage[pagewise,switch]{lineno} 

\theoremstyle{plain}
\newtheorem{theorem}{Theorem}[section]
\newtheorem{lemma}[theorem]{Lemma}

\newtheorem{corollary}[theorem]{Corollary}
\theoremstyle{definition}

\theoremstyle{remark}
\newtheorem{remark}[theorem]{Remark}

\newcommand{\Ha}{\operatorname{Ha}}

\title{On Semi-interval and Para-interval Graphs}
\author{
Chenxi Lv$^{1}$,
Yiming Lei$^{1}$\thanks{Corresponding author.
Email address: leiyiming008@gmail.com.},
Chunlei Xu$^{2}$\\[6pt]
\small $^{1}$College of Mathematical Sciences,
Bohai University, Jinzhou, China\\
\small $^{2}$College of Mathematics Science,
Inner Mongolia Minzu University, Tongliao, China
}
\date{}

\begin{document}
\maketitle

\begin{abstract}
Bau and Johnson introduced semi-interval graphs and asked whether one of the two extremal sets can always be replaced by a single vertex. A negative answer to this question is given by a tree which has no such representation. It is also shown that a finite connected graph with at least two vertices is a semi-interval graph if and only if it is bipartite. Para-interval graphs, for which one extremal set consists of a single vertex, are characterized by necessary and sufficient conditions on their distance matrices. A characterization of para-interval trees is obtained, and a relation between para-interval graphs and shortest-path interval graphs is established by identifying the vertices of an extremal set. Finally, leveled graphs and the corresponding Hasse diagrams of finite partially ordered sets are considered.
\end{abstract}



\section{Introduction}

Bau and Johnson~\cite{Bau} considered graphs which are unions of
shortest paths connecting two vertices. Such graphs were called
intervals. They also introduced semi-interval graphs by replacing
the two extremal vertices with two sets of vertices. In this paper,
the term interval graph is used in this shortest-path sense, rather
than in the usual sense of an intersection graph of intervals on
the real line.

Bau and Johnson proved that every connected bipartite graph is a
semi-interval graph~\cite[Lemma~3.2]{Bau}. They asked whether every
semi-interval graph has a representation in which one extremal set
consists of a single vertex. We give a counterexample to this
question. We then consider graphs which admit such a representation,
and call them para-interval graphs. We give a necessary and sufficient condition for
para-interval graphs in terms of distance matrices, and
apply this condition to characterize para-interval trees. We also
consider the relation of these graphs with intervals and finite
partially ordered sets. In the last two connections we shall use
the results and terminology of Bau and Johnson~\cite{Bau}.

\section{Semi-interval graphs}
\label{sec:semi}

Throughout this paper, graphs are finite and simple.
For basic concepts and terminology in graph theory,
the reader is referred to \cite{Bondy,Diestel}.

Let $G$ be a
connected graph with at least two vertices, and let $A,B\subseteq
V(G)$ be nonempty disjoint sets. Define
\[
\rho(A,B)=\min\{\rho(a,b):a\in A,\ b\in B\}.
\]
A shortest $A$--$B$ path is a path connecting a vertex of $A$ and
a vertex of $B$ whose length is $\rho(A,B)$. If there exist such
sets $A$ and $B$ such that $G$ is the union of all shortest
$A$--$B$ paths, then $G$ is called a semi-interval graph. The sets
$A$ and $B$ are called the extremal sets of the representation.
For $i\geq 0$, let
\[
L_i(A)=\{x\in V(G):\rho(A,x)=i\},\qquad
\rho(A,x)=\min_{a\in A}\rho(a,x).
\]

The sufficiency follows from Lemma~3.2 of Bau and
Johnson~\cite{Bau}.

\begin{theorem}\label{the:semi}
Let $G$ be a connected graph with at least two vertices. Then
$G$ is a semi-interval graph if and only if $G$ is bipartite.
\end{theorem}

\begin{proof}
Suppose first that $G$ is bipartite with bipartition $(A,B)$.
Then $\rho(A,B)=1$, and every edge of $G$ is a shortest
$A$--$B$ path. Hence $G$ is a semi-interval graph.

Suppose, conversely, that $G$ is a semi-interval graph with
extremal sets $A$ and $B$. Let $d=\rho(A,B)$ and let
$P=v_0v_1\cdots v_d$ be a shortest $A$--$B$ path, where
$v_0\in A$ and $v_d\in B$. For each $i$, we have
$\rho(A,v_i)\leq i$. If $\rho(A,v_i)<i$, then a path from $A$
to $v_i$ of length less than $i$, followed by
$v_iv_{i+1}\cdots v_d$, gives an $A$--$B$ walk of length less
than $d$. This is a contradiction. Hence $\rho(A,v_i)=i$.

Therefore, every edge of a shortest $A$--$B$ path joins two
consecutive distance levels. Since $G$ is the union of all such
paths, every edge of $G$ has this property. Let
\[
X=\bigcup_{i\text{ even}}L_i(A),\qquad
Y=\bigcup_{i\text{ odd}}L_i(A).
\]
Then $(X,Y)$ is a bipartition of $G$. Hence $G$ is bipartite.
\end{proof}

\begin{lemma}\label{lem:extremal-layers}
Let $G$ be a semi-interval graph with extremal sets $A$ and $B$.
If $d=\rho(A,B)$, then $B=L_d(A)$ and $A=L_d(B)$.
\end{lemma}

\begin{proof}
By the proof of Theorem~\ref{the:semi}, every shortest
$A$--$B$ path meets $L_i(A)$ at its $i$th vertex. Every
$b\in B$ belongs to such a path, and cannot be an internal
vertex of the path, since $\rho(A,b)\geq d$. Thus
$B\subseteq L_d(A)$. Conversely, every $x\in L_d(A)$ belongs
to a shortest $A$--$B$ path and can only be its last vertex.
Hence $L_d(A)\subseteq B$. Therefore, $B=L_d(A)$.
Interchanging $A$ and $B$ gives $A=L_d(B)$.
\end{proof}

\begin{figure}[htbp]
\centering
\begin{tikzpicture}[scale=1.05]
\coordinate (u) at (-0.65,0);
\coordinate (v) at (0.65,0);
\coordinate (a1) at (-1.50,0.78);
\coordinate (a2) at (-1.50,-0.78);
\coordinate (b1) at (1.50,0.78);
\coordinate (b2) at (1.50,-0.78);
\draw (a1)--(u)--(a2) (u)--(v) (v)--(b1) (v)--(b2);
\foreach \p in {u,v,a1,a2,b1,b2}{\filldraw[fill=white,draw=black] (\p) circle(2.1pt);}
\node[above] at (u) {$u$};
\node[above] at (v) {$v$};
\end{tikzpicture}
\caption{A semi-interval graph with no singleton extremal set.}
\label{fig:counterexample}
\end{figure}
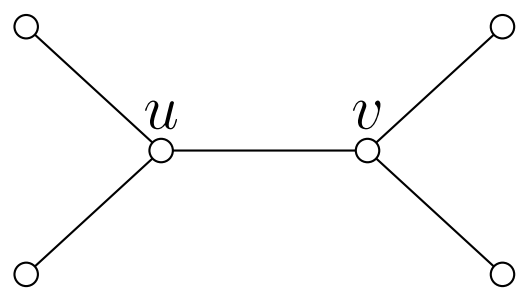

\begin{remark}\label{rem:counterexample}
Bau and Johnson~\cite{Bau} asked whether every semi-interval
graph has a representation with a singleton extremal set. We give
an example to show that this is not true.

Let $T$ be the tree in Figure~\ref{fig:counterexample}, consisting
of two adjacent vertices of degree $3$, each adjacent to two
leaves. Since $T$ is bipartite, it is a semi-interval graph by
Theorem~\ref{the:semi}.

Suppose that $T$ has a representation with extremal sets $\{s\}$
and $B$. Let $d=\rho(s,B)$. Every leaf different from $s$ must
be an end vertex of a shortest path connecting $s$ and $B$.
Hence every such leaf has distance $d$ from $s$.

If $s$ is a vertex of degree $3$, the four leaves are at
distances $1$ and $2$ from $s$. If $s$ is a leaf, the other
three leaves are at distances $2$ and $3$ from $s$. Thus, for
every choice of $s$, two leaves different from $s$ are at
different distances from $s$. This is a contradiction.
Therefore, $T$ has no representation with a singleton extremal
set. This gives a negative answer to the question of Bau and
Johnson.
\end{remark}

\section{Para-interval graphs}
\label{sec:para}

Let $G$ be a connected graph with at least two vertices. If
there exist $a\in V(G)$ and $B\subseteq V(G)$ such that $G$ is
the union of all shortest paths connecting $a$ and $B$, then
$G$ is called a para-interval graph. The vertex $a$ is called
the extremal vertex and $B$ the extremal set. Clearly, every
para-interval graph is a semi-interval graph.

For $s\in V(G)$ and each integer $i\geq 0$, let
$L_i(s)=\{x\in V(G):\rho(s,x)=i\}$.

\begin{lemma}\label{lem:ter1}
Let $G$ be a para-interval graph with extremal vertex $a$ and
extremal set $B$. If $d=\max\{\rho(a,x):x\in V(G)\}$, then
$B=L_d(a)$ and $d=\rho(a,B)$.
\end{lemma}

\begin{proof}
Let $r=\rho(a,B)$. By Lemma~\ref{lem:extremal-layers},
$B=L_r(a)$. Since every vertex of $G$ belongs to a shortest
path connecting $a$ and $B$, its distance from $a$ is at most
$r$. Hence $d\leq r$. As $B$ is nonempty and $B=L_r(a)$, we
also have $r\leq d$. Therefore, $d=r$ and $B=L_d(a)$.
\end{proof}

The preceding lemma shows that the extremal set is determined
by the extremal vertex. We now give a characterization of
para-interval graphs in terms of distance matrices.
Let $V(G)=\{1,2,\ldots,n\}$, and for each integer $k\geq 0$,
define $M_k=(m_{ij}^{(k)})_{n\times n}$ by
\[
m_{ij}^{(k)}=\begin{cases}
1,&\rho(i,j)=k,\\
0,&\rho(i,j)\neq k.
\end{cases}
\]
In particular, $M_1$ is the adjacency matrix of $G$.

\begin{theorem}\label{the:ter5}
Let $G$ be a connected graph with at least two vertices.
Then $G$ is a para-interval graph if and only if there
exists $a\in V(G)$, with
$d=\max\{\rho(a,x):x\in V(G)\}$, such that
\begin{equation}\label{eq:same-level}
(M_kM_1M_k)_{aa}=0\quad\text{for }0\leq k\leq d,
\end{equation}
and
\begin{equation}\label{eq:forward}
(M_{k+1}M_1-M_k)_{ax}\geq 0
\quad\text{for every }x\in V(G),\quad 0\leq k<d.
\end{equation}
Moreover, every vertex $a$ satisfying these conditions
is an extremal vertex of $G$.
\end{theorem}

\begin{proof}
Let $a\in V(G)$ and
$d=\max\{\rho(a,x):x\in V(G)\}$.
For $0\leq k\leq d$, by the definition of $M_k$,
\[
(M_kM_1M_k)_{aa}
=\sum_{x,y\in L_k(a)}m_{xy}^{(1)}
=2|E(G[L_k(a)])|.
\]
For $0\leq k<d$ and $x\in L_k(a)$, we also have
\[
(M_{k+1}M_1-M_k)_{ax}
=\sum_{y\in L_{k+1}(a)}m_{yx}^{(1)}-1
=|N(x)\cap L_{k+1}(a)|-1,
\]
where $N(x)$ denotes the set of vertices adjacent to $x$.
If $x\notin L_k(a)$, then
$(M_{k+1}M_1-M_k)_{ax}\geq 0$.

Suppose first that $G$ is a para-interval graph
with extremal vertex $a$ and extremal set $B$.
By Theorem~\ref{the:semi}, $G$ is bipartite.
Hence no two vertices of the same distance level
from $a$ are adjacent, and condition~\eqref{eq:same-level}
holds. By Lemma~\ref{lem:ter1}, $B=L_d(a)$.
For every $x\in L_k(a)$ with $k<d$, there exists
a shortest path connecting $a$ and some $b\in B$
which contains $x$. Hence $x$ has a neighbour in
$L_{k+1}(a)$, and condition~\eqref{eq:forward} holds.

Suppose, conversely, that there exists $a\in V(G)$
satisfying conditions~\eqref{eq:same-level} and
\eqref{eq:forward}. By condition~\eqref{eq:same-level},
no two vertices in the same distance level from $a$
are adjacent. Since adjacent vertices have distances
from $a$ which differ by at most one, every edge
of $G$ joins two consecutive distance levels.
By condition~\eqref{eq:forward}, every vertex in
$L_k(a)$, $0\leq k<d$, has a neighbour in
$L_{k+1}(a)$.

Let $B=L_d(a)$ and $xy\in E(G)$.
We may assume that $x\in L_i(a)$ and
$y\in L_{i+1}(a)$, where $0\leq i<d$.
By condition~\eqref{eq:forward}, there exists a path
of length $d-i-1$ connecting $y$ and some $b\in B$,
whose vertices belong successively to
$L_{i+1}(a),L_{i+2}(a),\ldots,L_d(a)$.
Since $x\in L_i(a)$, there exists a shortest path
of length $i$ connecting $a$ and $x$.
These two paths, together with the edge $xy$,
form a path of length $d$ connecting $a$ and $b$
which contains $xy$. Since $\rho(a,b)=d$,
this is a shortest path connecting $a$ and $b$.
Hence every edge of $G$ belongs to a shortest path
connecting $a$ and $B$. Therefore, $G$ is a
para-interval graph with extremal vertex $a$.
\end{proof}

\begin{theorem}\label{the:ter6}
Let $T$ be a tree with at least two vertices, and let
$a\in V(T)$. If $d=\max\{\rho(a,x):x\in V(T)\}$, then
$T$ is a para-interval graph with extremal vertex $a$ if
and only if every leaf of $T$ different from $a$ belongs
to $L_d(a)$.
\end{theorem}

\begin{proof}
Let $M_k$ be the distance matrices of $T$.
Suppose first that $T$ is a para-interval graph with
extremal vertex $a$. Let $x\neq a$ be a leaf. If
$x\in L_i(a)$ with $i<d$, then the only neighbour of
$x$ is its predecessor on the unique path connecting $a$
and $x$. Hence $(M_{i+1}M_1-M_i)_{ax}=-1$,
contrary to Theorem~\ref{the:ter5}.
Therefore, $x\in L_d(a)$.

Suppose, conversely, that every leaf different from $a$
belongs to $L_d(a)$. Let $x\in L_i(a)$ with $i<d$.
If $x=a$, then $x$ has a neighbour in $L_1(a)$.
Suppose that $x\neq a$. If $x$ has no neighbour in
$L_{i+1}(a)$, then its only neighbour is its predecessor
on the unique path connecting $a$ and $x$. Hence $x$
is a leaf, contrary to $i<d$. Therefore, every vertex
in $L_i(a)$, $0\leq i<d$, has a neighbour in
$L_{i+1}(a)$. Since no two vertices of the same
distance level are adjacent in a tree, both matrix
conditions of Theorem~\ref{the:ter5} hold.
Hence $T$ is a para-interval graph with extremal vertex $a$.
\end{proof}

Recall that a connected graph $G$ is called an
interval graph in the shortest-path sense if there
exist two vertices $a,b\in V(G)$ such that $G$ is
the union of all shortest paths connecting $a$ and
$b$ \cite{Bau}. The vertices $a$ and $b$ are called
the extremal vertices of the representation.
We now consider the relation between para-interval
graphs and interval graphs.

\begin{theorem}\label{the:ter7}
Let $G$ be a para-interval graph. Then there exists an
independent set $S\subseteq V(G)$ such that the graph
obtained by identifying all vertices of $S$ into a single
vertex, and suppressing multiple edges, is an interval
graph in the shortest-path sense.
\end{theorem}

\begin{proof}
Let $a$ and $B$ be the extremal vertex and extremal set
of $G$, respectively. Let $d=\max\{\rho(a,x):x\in V(G)\}$.
By Lemma~\ref{lem:ter1}, $B=L_d(a)$. Since $G$ is
bipartite, no two vertices of the same distance level are
adjacent. Hence $B$ is an independent set. Take $S=B$.

Identify all vertices of $B$ into a vertex $b^*$, and let
$G^*$ denote the resulting simple graph. Every neighbour
of a vertex in $B$ belongs to $L_{d-1}(a)$. Hence every
path from $a$ to $b^*$ in $G^*$ has length at least $d$.
On the other hand, a shortest path connecting $a$ to $B$
gives an $a$--$b^*$ path of length $d$. Therefore,
$\rho_{G^*}(a,b^*)=d$.

Every shortest path connecting $a$ to $B$ in $G$ gives a
shortest path connecting $a$ and $b^*$ in $G^*$. Since
$G$ is the union of all shortest paths connecting $a$
and $B$, every edge of $G^*$ belongs to a shortest
$a$--$b^*$ path. Hence $G^*$ is an interval graph.
\end{proof}

\begin{figure}[htbp]
\centering
\begin{tikzpicture}[scale=0.90]
\coordinate (u) at (0,0);
\coordinate (t) at (0,1);
\coordinate (v) at (0.85,-0.5);
\coordinate (w) at (0,-1);
\coordinate (x) at (0,-2);
\coordinate (y) at (0,-3);
\coordinate (a) at (-0.85,-0.22);
\coordinate (b) at (-0.85,-0.78);
\draw (t)--(u)--(v)--(w)--(x)--(y) (u)--(a) (w)--(b);
\foreach \p in {u,t,v,w,x,y,a,b}{\filldraw[fill=white,draw=black] (\p) circle(2.1pt);}
\node[left] at (a) {$a$};
\node[left] at (b) {$b$};
\node[right] at (v) {$c$};
\end{tikzpicture}
\caption{A graph which is not a para-interval graph.}
\label{fig:converse}
\end{figure}
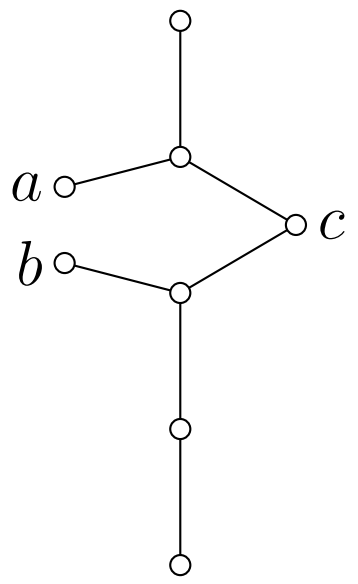

The converse of Theorem~\ref{the:ter7} is not true. Consider
the tree in Figure~\ref{fig:converse}. For any vertex $s$,
there are two leaves different from $s$ at different distances
from $s$. Hence by Theorem~\ref{the:ter6}, this tree is not a
para-interval graph. However, by identifying the independent
set $S=\{a,b,c\}$ into a single vertex, we obtain a path,
which is an interval graph.

As in Bau and Johnson~\cite{Bau}, a graph $G$ is called a
leveled graph if there exist $a,b\in V(G)$ such that every
vertex of $G$ belongs to some shortest path connecting $a$
and $b$. The following result is also a direct consequence of
Theorem~3.1 of Bau and Johnson~\cite{Bau}, which describes
the edges that may be added to an interval to obtain a leveled
graph.

\begin{theorem}\label{the:ter8}
Let $G$ be a leveled graph with vertices $a,b\in V(G)$.
Then $G$ is the union of all shortest $a$--$b$ paths if
and only if $\rho(a,u)\neq\rho(a,v)$ for every
$uv\in E(G)$.
\end{theorem}

\begin{proof}
Suppose that $\rho(a,u)\neq\rho(a,v)$ for every edge
$uv$. Since adjacent vertices have distances from $a$
which differ by at most one, every edge joins consecutive
distance levels. Let $uv\in E(G)$, and suppose that
$\rho(a,v)=\rho(a,u)+1=i+1$. Since $G$ is leveled,
$\rho(v,b)=\rho(a,b)-i-1$. A shortest path from $a$
to $u$, followed by $uv$ and a shortest path from $v$
to $b$, gives an $a$--$b$ walk of length $\rho(a,b)$.
Hence it is a shortest $a$--$b$ path containing $uv$.
Therefore, $G$ is the union of all shortest $a$--$b$
paths.

Suppose, conversely, that $G$ is the union of all
shortest $a$--$b$ paths. Along each such path, the
distance from $a$ increases by one at every step.
Hence no edge of $G$ joins two vertices at the same
distance from $a$.
\end{proof}

\section{Hasse diagrams}
\label{sec:hasse}

Let $\mathcal P=(V,\leq)$ be a finite partially ordered set.
We say that $y$ covers $x$ if $x<y$ and there is no
$z\in V$ such that $x<z<y$. Let $\Ha(\mathcal P)$ denote
the underlying undirected graph of the Hasse diagram of
$\mathcal P$. A maximal chain is a chain which cannot be
extended by adding another element. The length of a chain
is one less than its number of elements. We consider finite
posets whose maximal chains all have the same length.

Bau and Johnson~\cite{Bau} noted the relation between
shortest-path intervals and partially ordered sets with rank
functions, and considered the case of a unique minimal and a
unique maximal element in Corollary~4.5 of their paper.
Here we describe the partial order obtained from fixed
extremal sets of a semi-interval graph.

Let $G$ be a semi-interval graph with extremal sets $A$ and
$B$. Define $r(x)=\rho(A,x)$ for $x\in V(G)$, and define
a relation $\preceq$ on $V(G)$ by
\begin{equation}\label{eq:order}
x\preceq y\quad\Longleftrightarrow\quad
\rho(x,y)=r(y)-r(x).
\end{equation}

\begin{lemma}\label{lem:ha1}
The relation $\preceq$ defined by~\eqref{eq:order} is a
partial order on $V(G)$. Moreover, $y$ covers $x$ if and
only if $xy\in E(G)$ and $r(y)=r(x)+1$.
\end{lemma}

\begin{proof}
For every $x\in V(G)$, we have
$\rho(x,x)=0=r(x)-r(x)$, and hence $x\preceq x$.
If $x\preceq y$ and $y\preceq x$, then
$\rho(x,y)=r(y)-r(x)=r(x)-r(y)$, so that $x=y$.
Thus $\preceq$ is reflexive and antisymmetric.

Suppose that $x\preceq y$ and $y\preceq z$. By the
triangle inequality,
\[
r(z)-r(x)\leq\rho(x,z)\leq
\rho(x,y)+\rho(y,z)=r(z)-r(x).
\]
Hence $x\preceq z$, and the relation is transitive.
Therefore, it is a partial order.

Suppose that $y$ covers $x$. If $\rho(x,y)\geq 2$,
let $z$ be the neighbour of $x$ on a shortest path
connecting $x$ and $y$. Since
$\rho(x,y)=r(y)-r(x)$, we have
\[
r(y)\leq r(z)+\rho(z,y)
\leq r(x)+1+\rho(x,y)-1=r(y).
\]
Thus $r(z)=r(x)+1$ and
$\rho(z,y)=r(y)-r(z)$, giving $x\prec z\prec y$.
This is a contradiction. Hence $\rho(x,y)=1$ and
$r(y)=r(x)+1$.

Conversely, suppose that $xy\in E(G)$ and
$r(y)=r(x)+1$. Then $x\prec y$. There is no
$z$ with $x\prec z\prec y$, since this would give
$r(x)<r(z)<r(y)$. Hence $y$ covers $x$.
\end{proof}

\begin{theorem}\label{the:ha2}
Let $G$ be a connected graph with at least two vertices.
Then $G$ is a semi-interval graph if and only if there
exists a finite poset $\mathcal P$ such that all maximal
chains of $\mathcal P$ have the same length and $G$ is
isomorphic to $\Ha(\mathcal P)$. Moreover, if $G$ has
extremal sets $A$ and $B$, then $\mathcal P$ can be
chosen with minimal element set $A$ and maximal element
set $B$.
\end{theorem}

\begin{proof}
Suppose first that $G$ is a semi-interval graph with
extremal sets $A$ and $B$. Let $d=\rho(A,B)$ and
$r(x)=\rho(A,x)$. By Lemma~\ref{lem:ha1}, the relation
in~\eqref{eq:order} defines a finite poset $\mathcal P$.
By Lemma~\ref{lem:extremal-layers},
$A=L_0(A)$ and $B=L_d(A)$. Every vertex of $G$
belongs to a shortest $A$--$B$ path. Hence every
vertex outside $A$ has a neighbour in the preceding
level, and every vertex outside $B$ has a neighbour
in the succeeding level. By Lemma~\ref{lem:ha1}, the
minimal elements of $\mathcal P$ are exactly the
vertices of $A$, and the maximal elements are exactly
the vertices of $B$.

Every edge of $G$ joins consecutive distance levels.
By Lemma~\ref{lem:ha1}, these edges are precisely the
covering relations of $\mathcal P$. Hence $G$ is isomorphic to
$\Ha(\mathcal P)$. Let
$x_0\prec x_1\prec\cdots\prec x_t$ be a maximal chain.
Its first vertex is in $A$ and its last vertex is in
$B$. Since consecutive vertices form covering relations,
$r(x_{i+1})=r(x_i)+1$. Hence $t=d$, and every
maximal chain has length $d$.

Suppose, conversely, that $G$ is the underlying graph
of the Hasse diagram of a finite poset $\mathcal P$
whose maximal chains all have the same length $d$.
Let $A$ and $B$ be its minimal and maximal element
sets. Since $G$ is connected with at least two vertices,
$d\geq 1$ and $A\cap B=\varnothing$.

Every element $x$ has a well-defined rank $r(x)$,
namely the length of a saturated chain from a minimal
element to $x$. Indeed, if two such chains have lengths
$k$ and $k'$, extending both by the same saturated chain
from $x$ to a maximal element gives two maximal chains
of lengths $d$ and hence $k=k'$. Consequently, every
covering relation $x\prec y$ satisfies $r(y)=r(x)+1$.
In particular, $r(a)=0$ for $a\in A$ and $r(b)=d$ for
$b\in B$.

Every maximal chain determines an $A$--$B$ path of
length $d$. Since the rank changes by one along every
edge of $G$, every $A$--$B$ path has length at least
$d$. Therefore, $\rho(A,B)=d$, and every maximal chain
is a shortest $A$--$B$ path. Every covering relation
of a finite poset belongs to a maximal chain. Hence
every edge of $G$ belongs to a shortest $A$--$B$
path. Therefore, $G$ is a semi-interval graph.
\end{proof}

It follows from Theorem~\ref{the:ha2} that the condition
on the minimal and maximal element sets determines the
corresponding subclasses. The following statements also
recover the usual interpretation of a shortest-path interval
as a graded poset with a unique minimal and a unique maximal
element.

\begin{corollary}\label{cor:ha-para}
Let $G$ be a connected graph with at least two vertices.
Then $G$ is a para-interval graph if and only if there
exists a finite poset $\mathcal T$ with a unique minimal
element, with all maximal chains of the same length,
and with underlying Hasse graph isomorphic to $G$.
\end{corollary}

\begin{proof}
Suppose that $G$ is a para-interval graph with extremal
vertex $a$ and extremal set $B$. Applying
Theorem~\ref{the:ha2} to the extremal sets $\{a\}$ and
$B$ gives a poset $\mathcal T$ with unique minimal
element $a$. Conversely, if $\mathcal T$ has a unique
minimal element $a$, then, by Theorem~\ref{the:ha2},
$G$ is the union of shortest paths connecting $a$
and the maximal element set of $\mathcal T$. Hence
$G$ is a para-interval graph.
\end{proof}

\begin{corollary}\label{cor:ha-interval}
Let $G$ be a connected graph with at least two vertices.
Then $G$ is an interval graph in the shortest-path sense
if and only if there exists a finite poset $\mathcal Q$
with a unique minimal element, a unique maximal element,
and all maximal chains of the same length, such that
$G$ is isomorphic to the underlying Hasse graph of
$\mathcal Q$.
\end{corollary}

\begin{proof}
If $G$ is an interval with extremal vertices $a,b$,
then applying Theorem~\ref{the:ha2} with
$A=\{a\}$ and $B=\{b\}$ gives the required poset.
Conversely, if the minimal and maximal elements of
$\mathcal Q$ are $a$ and $b$, respectively, then
Theorem~\ref{the:ha2} shows that $G$ is the union
of all shortest $a$--$b$ paths. Hence $G$ is an
interval graph.
\end{proof}

\end{document}